\documentclass{article}

\usepackage{arxiv}

\usepackage[utf8]{inputenc} 
\usepackage[T1]{fontenc}    
\usepackage[hidelinks]{hyperref}       
\usepackage{url}            
\usepackage{booktabs}       
\usepackage{amsfonts}       
\usepackage{nicefrac}       
\usepackage{microtype}      
\usepackage{lipsum}		
\usepackage{graphicx}
\usepackage{doi}

\usepackage[backend=biber,style=numeric,maxnames=4,minnames=4,url=false,isbn=false,sorting=none]{biblatex}
\usepackage{amsmath,amssymb}
\usepackage{multirow}
\usepackage{mathtools}
\usepackage{colortbl}
\usepackage{stmaryrd}
\usepackage{pgfplotstable}
\usepackage{booktabs}
\usepackage{cleveref}
\pgfplotsset{compat=1.18}
\crefname{figure}{figure}{figures}  
\usetikzlibrary{calc}

\usepackage[acronym]{glossaries}
\newacronym[shortplural={BCs},longplural={Boundary Conditions}]{bc}{BC}{Boundary Condition}
\newacronym[shortplural={DOFs},longplural={degrees of freedom}]{dof}{DOF}{degree of freedom}
\newacronym{wss}{WSS}{wall shear stress}
\newacronym{fe}{FE}{finite element}
\newacronym{fem}{FEM}{Finite Element Method}
\newacronym{cg1}{P1}{First order Lagrange}
\newacronym{dg1}{DG-1}{Discontinuous Galerkin first order}
\newacronym{dg0}{DG-0}{Discontinuous Galerkin zero order}
\newacronym{lsa}{LSA}{low shear area}
\newacronym{cfd}{CFD}{computational fluid dynamics}

\providecommand{\ix}[1]{_{\textrm{#1}}} 
\providecommand{\na}[1]{^{\textrm{#1}}} 
\providecommand{\dx}{\hspace{2pt}\text{dx}\hspace{2pt}} 
\providecommand{\ds}{\hspace{2pt}\text{dS}\hspace{2pt}} 
\DeclareMathOperator{\divergenceOperator}{div}
\renewcommand{\vec}[1]{\ensuremath{\boldsymbol{#1}}}
\newcommand\defeq{\mathrel{\stackrel{\makebox[0pt]{\mbox{\normalfont\tiny def}}}{=}}}
\newcommand{\wss}{\boldsymbol{\tau}}
\newcommand{\wssb}{\boldsymbol{\tau_b}}

\title{On the numerical evaluation of wall shear stress using the finite element method}

\author{ \href{https://orcid.org/0009-0000-5948-233X}{\includegraphics[scale=0.06]{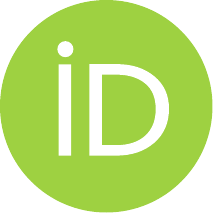}\hspace{1mm}Jana Brun\'atov\'a$^{1,2}$}\\
	$^1$Mathematical Institute, Charles University\\
	$^2$Bernoulli Institute, University of Groningen \\
	\texttt{brunatova@karlin.mff.cuni.cz} \\
    \And
	\href{https://orcid.org/0000-0001-6489-8858}{\includegraphics[scale=0.06]{orcid.pdf}\hspace{1mm}Jørgen Schartum Dokken$^3$} \\
    $^3$Scientific Computing and Numerical Analysis,\\
	Simula Research Laboratory\\
	\texttt{dokken@simula.no} \\
    \And
	\href{https://orcid.org/0000-0002-2907-0171}{\includegraphics[scale=0.06]{orcid.pdf}\hspace{1mm}Kristian Valen-Sendstad$^4$} \\
    $^4$Computational Physiology,\\
	Simula Research Laboratory\\
	\texttt{kvs@simula.no} \\
	\And
	\href{https://orcid.org/0000-0001-5862-2353}{\includegraphics[scale=0.06]{orcid.pdf}\hspace{1mm}Jaroslav Hron$^1$} \\
	$^1$Mathematical Institute, Charles University\\
	\texttt{hron@karlin.mff.cuni.cz} \\
}

\renewcommand{\shorttitle}{On the numerical evaluation of wall shear stress using the finite element method}

\hypersetup{
pdftitle={On the numerical evaluation of wall shear stress using the finite element method},
pdfsubject={math.NA},
pdfauthor={Jana~Brun\'atov\'a, Jørgen Schartum Dokken, Kristian Valen-Sendstad, Jaroslav Hron},
pdfkeywords={Wall shear stress, Finite element method, Boundary-flux evaluation,Stokes flow, Navier--Stokes equations},
}

\begin{document}
\maketitle

\begin{abstract}
\Gls{wss} is a crucial hemodynamic quantity extensively studied in cardiovascular research, yet its numerical computation is not straightforward. 
This work aims to compare \gls{wss} results obtained from two different finite element discretizations, quantify the differences between continuous and discontinuous stresses, and introduce a novel method for \gls{wss} evaluation through the formulation of a boundary-flux problem.

Two benchmark problems are considered -- a 2D Stokes flow on a unit square and a 3D Poiseuille flow through a cylindrical pipe.
These are followed by investigations of steady-state Navier--Stokes flow in two patient-specific aneurysms.
The study focuses on P1/P1 stabilized and Taylor–Hood P2/P1 mixed finite elements for velocity and pressure.
\gls{wss} is computed using either the proposed boundary-flux method or as a projection of tangential traction onto \gls{cg1}, \gls{dg1}, or \gls{dg0} space.

For the P1/P1 stabilized element, the boundary-flux and \gls{cg1} projection methods yielded equivalent results. 
With the P2/P1 element, the boundary-flux evaluation demonstrated faster convergence in the Poiseuille flow example but showed increased sensitivity to pressure field inaccuracies in patient-specific geometries compared to the projection method. 
In patient-specific cases, the P2/P1 element exhibited superior robustness to mesh size when evaluating average \gls{wss} and \gls{lsa}, outperforming the P1/P1 stabilized element.

Projecting discontinuous finite element results into continuous spaces can introduce artifacts, such as the Gibbs phenomenon. 
Consequently, it becomes crucial to carefully select the finite element space for boundary stress calculations -- not only in applications involving \gls{wss} computations for aneurysms.
\end{abstract}


\section{Introduction}\label{introduction}
 
Blood flow-induced \acrfull{wss}, the tangential component of stress acting on the endothelium, is essential for maintaining endothelial cell function under normal physiological conditions~\cite{chien2007mechanotransduction}.
Abnormal \gls{wss} patterns have been associated with the initiation, progression, and outcomes of vascular diseases such as atherosclerosis~\cite{malek1999hemodynamic} and aneurysms~\cite{Frosen2019}.
Over the past two decades, medical-image based \gls{cfd} has played an important role in estimating these stresses~\cite{steinman2003image}. However, there are conflicting results  about whether physiologically too high~\cite{cebral2011quantitative} or too low~\cite{Xiang2011} \gls{wss} is associated with aneurysm rupture status in large retrospective studies. A hypothesis has been proposed
to unify these apparently conflicting results~\cite{Meng2013highlowWSS}, but one of the fundamental problems is that the \gls{wss} computations and their mathematical definitions vary considerably~\cite{liang2019towards}.

Although many studies have investigated \gls{wss} in aneurysms -- see recent review articles \cite{Zhou2017review, Wang2024} -- a consensus on the optimal numerical approach for its evaluation has yet to be established.
The lack of standardization contributes to variability in patient-specific CFD results, as reported in~\cite{ValenSendstad2018}, which undermines the trustworthiness required for eventual clinical translation.
Consequently, developing a reliable and accurate image-based CFD pipeline for performing patient-specific simulations is crucial.
Advancements in this area would enhance the reproducibility of computational studies, facilitate comparisons between studies and aid in translating computational findings into clinical practice.
This study concentrates on the post-processing phase of the CFD pipeline.

The main objective of this study is to examine numerical methods for \gls{wss} assessment within patient-specific \gls{cfd} models.
In addition to the standard method for \gls{wss} computation, we propose and investigate an alternative approach based on the formulation of a boundary flux problem similar to \cite{Carey1985,Brummelen2011}.
Furthermore, this study investigates the impact of the choice of finite element spaces on the accuracy of \gls{wss} assessment, which has not been shown in previous studies.

All numerical methods are validated against two benchmark problems with analytical solutions: 2D Stokes flow on a unit square and 3D Poiseuille flow in a cylindrical pipe.
Subsequently, the same approach is applied to two patient-specific aneurysm geometries to demonstrate its ability to handle anatomically relevant situations under simplified flow conditions.
Numerical simulations of both Stokes and Navier--Stokes flows employ a stable mixed finite element space: either the P1/P1 element with appropriate stabilization or the Taylor–Hood P2/P1 element.
Our open-source code based on FEniCS library~\cite{alnaes2015fenics} is used for numerical simulations.

The findings of this study aim to elucidate the numerical aspects of \gls{wss} assessment, providing a comparative evaluation of the methods employed and discussing the applicability of each approach. 
By examining the differences among these techniques, we aim to offer insights into their suitability for accurate and reliable \gls{wss} evaluation in patient-specific \gls{cfd} modeling.

\section{Theory}\label{theory}

In this study, we consider incompressible fluid flow, modeled by either the Stokes equation or the Navier--Stokes equations.
The \gls{fem} will be used for spatial discretization, and a finite difference scheme for temporal discretization.
\gls{wss} will be assessed in two ways; a projection of the tangential traction force into an appropriate finite element space, and a boundary-flux evaluation method.
Detailed descriptions of both assessment techniques are provided in this section.
In the following, we will denote the domain occupied by the fluid as $\Omega\subset\mathbb{R}^n$, where $n=2,3$.

\subsection{Stationary Stokes flow}
The stationary Stokes flow of an incompressible fluid is governed by the following equations
\begin{subequations}
  \begin{align}
   \hspace{4.0cm} 
   \divergenceOperator \mathbb{T}(\vec v, p) + \vec{f}&= \mathbf{0} &\text{in }& \Omega, \\
    \mathbb{T}(\vec v, p) &= -p \mathbb{I} + \nu \nabla \vec{v} &\text{in }& \Omega, \label{eqn:Cauchy_stress_Stokes}\\
   \divergenceOperator \vec{v} 
    &=
    0 &\text{in }& \Omega,
   \end{align}
   \label{eqn:Stokes_eqns}
\end{subequations}
where $\vec v$ is the velocity and $\nu$ kinematic viscosity of the fluid; $p$ denotes the pressure; $\mathbb{T}(\vec v, p) = -p \mathbb{I} + \nu \nabla \vec{v}$ is the Cauchy stress tensor, where $\mathbb I$ is the identity matrix in dimension 2 or 3; and $\vec f$ denotes any external volume forces.
\glspl{bc} are usually imposed by setting the velocity on the boundary (Dirichlet \gls{bc}), by prescribing traction force acting on the boundary (Neumann \gls{bc}), or by a combination of these two
\begin{equation}
\text{\glspl{bc}}
\begin{cases}
    \vec{v} = \vec{v\ix{bc}}, &\text{on } \Gamma\ix D,\\    
    \mathbb{T}(\vec v, p) \vec{n} = \vec{g} &\text{on } \Gamma \ix N,
\end{cases}
\label{eqn:BCs_generally}
\end{equation}
where $\overline{\Gamma\ix D \cup \Gamma\ix N} = \partial \Omega$, $\Gamma\ix D \cap \Gamma\ix N = \varnothing$.
When imposing a Dirichlet \gls{bc} for velocity, it can be applied either in a strong sense~\cite{Chouly2024bcs}, which is the common practice, or in a weak sense using Nitsche's method~\cite{Nitsche1971,Burman2012_Nitsche,blank2018stokes}.
The choice of whether to impose the \gls{bc} strongly or weakly will depend on the selection of finite element for velocity and pressure.
Additionally, some stabilization terms may be included, depending on the stability requirements of the finite element formulation.
The corresponding weak formulation of the problem \ref{eqn:Stokes_eqns} follows
\begin{align}
\begin{split}
    \int_\Omega \mathbb{T} : \nabla\vec{w} \dx &- \int_{\Gamma\ix N} \vec g \cdot \vec{w} \ds + \int_\Omega q \divergenceOperator \vec{v} \dx\\
    &+ \mathcal{N}_{\Gamma\ix D}[(\vec{v}\ix{bc};\vec{v},p),(\vec{w},q)] 
    + \mathcal{S}[(\vec{v},p),(\vec{w},q)]
    = \int_\Omega \vec f \cdot \vec w \dx \qquad \forall (\vec{w}, q) \in \textbf{V}\times Q,
\end{split}\label{eqn:Stokes_weak_form}
\end{align}
where we used integration by parts and where $\mathcal{N}_{\Gamma\ix D}$ denotes additional terms if \glspl{bc} on $\Gamma\ix D$ are imposed by Nitsche's method, while $\mathcal{S}$ denotes stabilization terms, which will be described in Section \ref{methods}. 
In the case of Dirichlet \glspl{bc} only, the function space for pressure is defined as $Q \defeq L_0^2(\Omega)$, representing $L^2$ functions with zero integral over the domain $\Omega$.
When a Neumann \gls{bc} is included, which will be the case for 3D problems, the function space becomes $Q \defeq L^2(\Omega)$.
The function space for velocity will be denoted differently when using Nitsche enforcement of \gls{bc}, $\vec{V} \defeq \vec{W}^{1,2}(\Omega)$, and strong enforcement of \gls{bc}, $\vec{V}^{\vec{v}\ix{bc}} \defeq \{\vec{v} \in \vec{W}^{1,2}(\Omega)~\vert~ \vec v = \vec v\ix{bc} \text{ on } \Gamma\ix D\}$.

We will utilize the non-symmetric Nitsche's method to enforce the Dirichlet \glspl{bc}\cite{blank2018stokes}.
In principle, the non-symmetric Nitsche's method does not need a penalization term like its symmetric variant; however, it has been shown that for Navier slip \gls{bc}, using the penalization term can significantly improve the accuracy of the quantities of interest~\cite{Chabiniok2021}. 
The non-symmetric Nitsche method with penalization parameter for enforcing Dirichlet \glspl{bc} on boundary $\Gamma\ix D$ follows
\begin{gather}
    \mathcal{N}_{\Gamma\ix D}[(\vec{v}\ix{bc};\vec{v},p),(\vec{w},q)] = 
    - \int_{\Gamma\ix{D}} [\mathbb{T}(\vec{v},p)\vec{n}]\cdot\vec{w} \ds 
    + \int_{\Gamma\ix{D}} [\mathbb{T}(\vec{w},q)\vec{n}]\cdot (\vec{v} - \vec{v}\ix{bc}) \ds 
    + \frac{\beta \nu}{h}\int_{\Gamma\ix{D}} (\vec{v} - \vec{v}\ix{bc})\cdot\vec{w} \ds,
\end{gather}
where $h$ denotes the circumdiameter, i.e. the diameter of a circumscribed sphere for each cell, and $\beta$ is a parameter that needs to be adjusted for the particular problem.

The finite element choice is crucial in obtaining stable numerical solutions. Selecting an inf-sup stable finite element pair is necessary~\cite{Boffi2013}, or, alternatively, implementing stabilization techniques when the pair lacks inherent inf-sup stability.
In this study, we focus on two finite elements: the inherently inf-sup stable Taylor–Hood P2/P1 element and the P1/P1 stabilized element.

\subsection{Navier--Stokes flow}

The blood flow within a vessel is governed by the incompressible Navier--Stokes equations under the usual assumptions that the vessel walls are rigid and impermeable, and that volume forces, such as gravitational effects, are neglected. 
The problem follows
\begin{subequations}
  \begin{align}
   \hspace{3.0cm} \rho \frac{\partial \vec{v}}{\partial t} 
    + \rho \left( \nabla \vec{v} \right) \vec{v}
    &= \divergenceOperator \mathbb{T}(\vec v, p) &\text{in }& \Omega \ix{T}, \label{eqn:general_N-S} \\
    \mathbb{T}(\vec v, p)
    &= -p \mathbb{I} + 2 \mu \mathbb{D} &\text{in }& \Omega \ix{T}, \label{eqn:Cauchy_stress}\\
   \divergenceOperator \vec{v} 
    &=
    0 &\text{in }& \Omega \ix{T},\label{eqn:div_free}\\
    \vec{v}(0,x)
    &=
    \vec{v}_0 &\text{in }& \Omega,\label{eqn:ic}\\
    \vec{v}
    &= \vec{v_{in}}(t, x), &\text{on }& \Gamma\ix{in}\na T,\\
    \vec{v}
    &= \vec{0} &\text{on }& \Gamma\ix{wall}\na T,\label{eqn:no_slilp}\\
    \mathbb{T} \vec{n}   
    &= \vec{0} &\text{on }& \Gamma\ix{out}\na T, \label{eqn:do_nothing}
   \end{align}
   \label{eqn:problem}
\end{subequations}
where $\rho$ is the density and $\mu$ the dynamic viscosity of the fluid; $t$ denotes time; $\mathbb T$ denotes the Cauchy stress tensor, where $\mathbb I$ is the identity matrix in dimension 3, and $\mathbb D = \frac{1}{2} \bigl(\nabla \vec v + (\nabla \vec v)\na T\bigr)$ is the symmetric part of the velocity gradient.
The initial condition $\vec{v}_0$ is usually set to zero velocity everywhere; the inflow velocity is prescribed by a function $\vec v \ix{in}$ depending on time and spatial coordinates; the no-slip \gls{bc} is set on walls, and the do-nothing \gls{bc} on all outlets.
The computational domain $\Omega\subset\mathbb{R}^3$ in time is denoted $\Omega \ix{T} \defeq (0,T) \times \Omega$, since the simulation starts at $t=0$ and ends at $t=T$. Similarly, we denote the inlet surface by $\Gamma\ix{in}$, outlet parts by $\Gamma\ix{out}$ and the walls by $\Gamma\ix{wall}$.

The approximate solution of the problem \eqref{eqn:problem} must be sought in a weak sense and discretized in both space and time. Function spaces for velocity and pressure are defined similarly as in the Stokes example: 
$\vec{V}^{\vec{v}\ix{bc}} \defeq \{\vec{v} \in \textbf{W}^{1,2}(\Omega)~\vert~ \vec v = \vec v\ix{in}\text{ on } \Gamma\ix{in}\na T \wedge \vec v = \vec 0\text{ on }\Gamma\ix{wall}\na T\}$ using strong enforcement of Dirichlet \gls{bc}, while
$\textbf{V} \defeq \{\vec{v} \in \textbf{W}^{1,2}(\Omega)\}$ for Nitsche enforcement; and $Q \defeq L^2(\Omega)$.
By employing integration by parts, the corresponding weak formulation yields
\begin{align}
\begin{split}
    &\rho \int_{\Omega} \left(\frac{\partial \vec{v}}{\partial t}  
    + \left( \nabla \vec{v} \right) \vec{v}\right) \cdot \vec{w} \dx 
    + \int_{\Omega} \mathbb{T}:\nabla \vec{w}\dx 
    + \int_{\Omega} q \divergenceOperator \vec{v} \dx
    + \mathcal{S}[(\vec{v},p),(\vec{w},q)]\\
    &+ \mathcal{N}_{\Gamma\ix{in}}[(\vec{v}\ix{in};\vec{v},p),(\vec{w},q)]
    + \mathcal{N}_{\Gamma\ix{wall}}[(\vec{0};\vec{v},p),(\vec{w},q)]
    = 0 \qquad 
    \text{a.e. in }(0,T) \quad
    \forall (\vec{w}, q) \in \textbf{V}\times Q,
\end{split}\label{eqn:weak_form}
\end{align}
where $\mathcal{N}_{\Gamma\ix{in}}$ and $\mathcal{N}_{\Gamma\ix{wall}}$ denotes Nitsche terms for enforcing the inflow and wall \glspl{bc}, respectively. The last term $\mathcal{S}$ includes possible stabilization terms which are particularly crucial when the chosen finite element pair for velocity and pressure does not inherently ensure stability.

\subsection{Wall shear stress evaluation}
\label{sec:wss_evaluation}

\gls{wss} is defined as the tangential component of the traction force per unit area exerted by a flowing fluid on a boundary. 
In the SI unit system, \gls{wss} is measured in pascals (Pa). 
However, it is most commonly treated as a vectorial quantity, representing both the magnitude and direction of the tangential force. By introducing the notation for tangential part of a vector $\vec a_t = \vec a - (\vec a \cdot \vec n) \vec n$, the definition of \gls{wss}, which will be denoted $\wss$ hereafter, follows
\begin{equation}
    \vec{\wss} \defeq
    (\mathbb{T}\vec n)_t 
    = \bigl((- p \mathbb{I} + 2 \mu \mathbb{D})\vec n \bigr)_t.
\end{equation}
Some definitions of \gls{wss} in the literature exclude the pressure term, as this component vanishes when the tangential projection is applied, as shown in Equation~\eqref{eq:wssproc}.

\subsubsection{WSS computation by projection}
\label{sec:projection_method}

Following the definition, \gls{wss} is typically numerically evaluated by performing an $L_2$ projection to an appropriate function space $\boldsymbol{S}$
\begin{equation}\label{eq:wssproc}
    \int_{\Gamma\ix{wall}} \wss \cdot \boldsymbol{\phi} \ds = 
    \int_{\Gamma\ix{wall}} \bigl(\mathbb{T}\vec n - (\mathbb{T}\vec n \cdot \vec n) \vec n \bigr)\cdot \boldsymbol{\phi} \ds = 2 \mu \int_{\Gamma\ix{wall}}  \bigl((\mathbb{D}\vec n) - (\mathbb{D}\vec n \cdot \vec n) \vec n\bigr) \cdot \boldsymbol{\phi} \ds \quad \forall \boldsymbol{\phi} \in \boldsymbol{S}.
\end{equation}

The function space for the above test functions, as well as for \gls{wss}, is often set to vectorial \gls{cg1} space. For example, this approach is implemented in the open-source VaMPy package \cite{Kjeldsberg2023} and the SimVascular software \cite{Updegrove2016}.
Commercial software typically does not specify the computational methods used; however, output data appear to be continuous and linear, as indicated by studies utilizing COMSOL Multiphysics \cite{Lee2018} and STAR-CCM+ \cite{Febina2018}, both of which are based on \gls{fem}. 
Although the \gls{cg1} space seems like an intuitively good option, the continuity of the velocity field across elements does not ensure the continuity of velocity gradients or traction forces.
Therefore, we hypothesize that utilizing discontinuous Galerkin polynomial function spaces may offer benefits for \gls{wss} computation.
Depending on the choice of \gls{fe} space for velocities, the natural discontinuous space for \gls{wss} should be of one order less than for velocity.
For instance, \gls{dg0} \gls{wss} when using P1/P1 stabilized element, while \gls{dg1} for P2/P1 element.
Additionally, a significant benefit of this approach is that the projection mass matrix can be inverted locally, as there is no coupling across elements. This local inversion can substantially accelerate computations on fine meshes, which might be beneficial especially in the context of moving domains where the mass matrix must be recomputed with each mesh movement.

Moreover, it is important whether the projection is done only on the boundaries of the domain $\boldsymbol{S}= \boldsymbol{S}(\Gamma\ix{wall})$ or in the entire domain $\boldsymbol{S}=\boldsymbol{S}(\Omega)$. 
In the latter case, one commonly sets \gls{wss} to zero at all interior \glspl{dof}. This adjustment does not follow from the formulation above, as those \glspl{dof} are not part of the equation. Without this notion, the corresponding matrix A would be singular.
It is clear that the test functions will be different depending on the space $\boldsymbol{S}$, which can especially influence results in corner elements.

\subsubsection{Boundary-flux evaluation of WSS}
\label{sec:boundary-flux_evaluation}

The boundary-flux evaluation technique was first introduced by Carey in 1985 for the Poisson problem~\cite{Carey1985} and was later extended by Brummelen in 2011 to address boundary-coupled problems~\cite{Brummelen2011}.
This technique enables the calculation of fluxes or stresses directly from finite element solutions.
It is an essential technique, for example, in fluid-structure interaction problems where the flux evaluation appears implicitly in the problem formulation and it has a direct connection to methods using Lagrange multipliers to enforce Dirichlet boundary conditions~\cite{Babuska1973,Bramble1981}.
In our case, the stress on the boundary is reconstructed from a solution to the Navier-Stokes equations as a postprocessing operation.

Assume that a weak solution $(\vec{v}, p)$ to Equation~\eqref{eqn:Stokes_weak_form} has been obtained under Dirichlet \glspl{bc} applied to all boundaries, resulting in $\Gamma_N = \varnothing$. The objective is to determine the flux through the walls, which upholds the no-slip boundary condition. This flux is identified with the traction force on the boundary, denoted by $\vec{t} = \mathbb{T}(\vec{v}, p)\vec{n}$. Mathematically, this approach is equivalent to prescribing Neumann \gls{bc} in place of Dirichlet \gls{bc} on all boundaries and subsequently evaluating the force on the boundary using the known solution $(\vec{v}, p)$; see\cite{Carey1985, Brummelen2011} for further details.
Moreover, let us assume that the Dirichlet \gls{bc} was imposed strongly in Equation \eqref{eqn:Stokes_weak_form}.
The term $\mathcal{N}_{\partial\Omega}[(\vec{v}\ix{bc};\vec{v},p),(\vec{w},q)]$ therefore vanishes and the boundary-flux evaluation for traction force reads
\begin{gather}
\label{eqn:boundary_flux_DirichletBC}
    \int_{\partial \Omega} \vec t \cdot \boldsymbol{\phi} \ds = 
    \int_{\Omega} \mathbb{T}(\vec v, p):\nabla \boldsymbol{\phi}\dx 
    - \int_\Omega \vec f \cdot \boldsymbol{\phi} \dx 
    + \mathcal{S}[(\vec{v},p),(\boldsymbol{\phi},0)]
    \qquad \forall \boldsymbol{\phi} \in \textbf{V}.
\end{gather}
It should be noted that the term associated with the test function from pressure space, $q$, is set to zero due to the absence of a corresponding term on the left-hand side.
Furthermore, the objective is to compute \gls{wss}, which is solely the tangential component of $\vec t$ that will be denoted by $\wssb$ below. 
By subtracting the normal traction from both sides of Equation \eqref{eqn:boundary_flux_DirichletBC}, one obtains the boundary-flux evaluation for \gls{wss} as follows
\begin{gather}
    \int_{\partial\Omega} \wssb \cdot \boldsymbol{\phi} \ds = 
    \int_{\Omega} \mathbb{T}(\vec v, p):\nabla \boldsymbol{\phi}\dx
    - \int_\Omega \vec f \cdot \boldsymbol{\phi} \dx 
    + \mathcal{S}[(\vec{v},p),(\boldsymbol{\phi},0)]
    - \int_{\partial\Omega} \bigl(\mathbb{T}(\vec v, p)\vec{n} \cdot \vec n \bigr) \vec n \cdot \boldsymbol{\phi} \ds
     \qquad \forall \boldsymbol{\phi} \in \textbf{V}.
\end{gather}

Subsequently, we turn our attention to the problem \eqref{eqn:problem}, with the aim of evaluating \gls{wss} only on walls ($\Gamma\ix{wall}$).
Consequently, the right-hand side must include both the time derivative term and the convective term. 
Moreover, the following Neumann boundary terms corresponding to the remaining boundary segments, $\Gamma\ix{in}$ and $\Gamma\ix{out}$, must be incorporated into the right-hand side of the boundary-flux evaluation, regardless of the original \gls{bc} type
\begin{align}\label{eqn:wss_weak_Dirichlet}
\begin{split}
    \int_{\Gamma\ix{wall}} \wssb \cdot \boldsymbol{\phi} \ds =&
    \rho \int_{\Omega} \left(\frac{\partial \vec{v}}{\partial t}  
    + \left( \nabla \vec{v} \right) \vec{v}\right) \cdot \boldsymbol{\phi}\dx 
    + \int_{\Omega} \mathbb{T}(\vec v, p):\nabla \boldsymbol{\phi}\dx
    - \int_\Omega \vec f \cdot \boldsymbol{\phi} \dx 
    + \mathcal{S}[(\vec{v},p),(\boldsymbol{\phi},0)] \\
    &- \int_{\Gamma\ix{wall}} \bigl(\mathbb{T}(\vec v, p)\vec{n} \cdot \vec n \bigr) \vec n \cdot \boldsymbol{\phi} \ds
    - \int_{\Gamma\ix{in}} \mathbb{T}(\vec v, p)\vec{n} \cdot \boldsymbol{\phi} \ds     
    - \int_{\Gamma\ix{out}} \mathbb{T}(\vec v, p)\vec{n} \cdot \boldsymbol{\phi} \ds   \qquad \forall \text{$\boldsymbol{\phi} \in \textbf{V}$.}
\end{split}
\end{align}
Note that the natural outlet \gls{bc} as well as the source term $\vec f$ might be omitted thanks to our original problem formulation, and the time derivative term may be omitted as well in the case of steady-state flow.
However, we keep all terms in the formulation for the sake of completeness.

A generalization of the above boundary-flux evaluation is necessary when the Dirichlet \gls{bc} on the wall was imposed by Nitsche's method.
In that case, two of the Nitsche terms remain on the right-hand side, while the third term is the unknown flux.
The boundary-flux evaluation for \gls{wss} thus reads
\begin{equation}\label{def:wss_weak}
\boxed{
\begin{aligned}
\mathclap{\text{Find $\wssb \in \boldsymbol{L}^2(\Gamma\ix{wall})$ such that:}}\\
    \int_{\Gamma\ix{wall}} \wssb \cdot \boldsymbol{\phi} \ds =& 
    \rho \int_{\Omega} \left(\frac{\partial \vec{v}}{\partial t}  
    + \left( \nabla \vec{v} \right) \vec{v}\right) \cdot \boldsymbol{\phi} \dx 
    + \int_{\Omega} \mathbb{T}(\vec v, p):\nabla \boldsymbol{\phi}\dx
    - \int_\Omega \vec f \cdot \boldsymbol{\phi} \dx 
    + \mathcal{S}[(\vec{v},p),(\boldsymbol{\phi},0)] \\
    &- \int_{\Gamma\ix{wall}} \bigl(\mathbb{T}(\vec v, p)\vec{n} \cdot \vec n \bigr) \vec n \cdot \boldsymbol{\phi} \ds
    + \int_{\Gamma\ix{wall}} \bigl(\mathbb{T}(\boldsymbol{\phi},0)\vec{n}\bigr)\cdot (\vec{v} - \vec{v}\ix{bc}) \ds 
    + \frac{\beta \mu}{h}\int_{\Gamma\ix{wall}} (\vec{v} - \vec{v}\ix{bc})\cdot\boldsymbol{\phi} \ds\\
    &- \int_{\Gamma\ix{in}} \mathbb{T}(\vec v, p)\vec{n} \cdot \boldsymbol{\phi} \ds
    - \int_{\Gamma\ix{out}} \mathbb{T}(\vec v, p)\vec{n} \cdot \boldsymbol{\phi} \ds
    \qquad \forall \text{$\boldsymbol{\phi} \in \textbf{V}$.}
\end{aligned}
}
\end{equation}

\section{Methods}\label{methods}
In this study, two different techniques for \gls{wss} evaluation are examined -- the projection method and the boundary-flux evaluation method described in Section \ref{sec:wss_evaluation}. 
For the projection method, three different finite element spaces for \gls{wss} are considered: the \gls{dg0}, \gls{dg1} and \gls{cg1} spaces. For the boundary-flux evaluation, we consider \gls{cg1} space for both finite elements.
As mentioned above, two types of finite elements for velocity and pressure are investigated: the Taylor--Hood P2/P1 element and stabilized P1/P1 element.
Taylor--Hood P2/P1 element satisfies the inf-sup condition~\cite{Boffi2013} and therefore does not require additional stabilization, however, the computational cost is demanding.
In contrast, the P1/P1 element is computationally cheaper but it violates the inf-sup condition and therefore needs to be stabilized; the specific stabilization method depends on the problem at hand, as described below. 

The accuracy of the aforementioned methods and finite element choices is assessed using two academic examples with known solutions and one real-world example.
The first example involves stationary Stokes flow in a 2D unit square domain, and the second examines 3D Poiseuille flow through a cylindrical pipe.
In these cases, the exact solutions are known, allowing for the assessment of convergence rates by measuring the relative errors with respect to the exact solutions.
The final example focuses on stationary Navier--Stokes flow in 3D patient-specific aneurysm geometries.
Since exact solutions are not available for these scenarios, the maximum, minimum and average \gls{wss} values over the aneurysm dome, together with the \gls{lsa} indicator, are shown for various mesh sizes.

\subsection{Implementation}
All computational codes used in this study are publicly available on Zenodo \href{https://doi.org/10.5281/zenodo.14506052}{doi.org/10.5281/zenodo.14506052}.
Additionally, all computational meshes can be found at \href{https://doi.org/10.5281/zenodo.14503385}{doi.org/10.5281/zenodo.14503385}, which ensures the reproducibility of our results.
A monolithic scheme is utilized to obtain numerical solutions of Stokes and Navier--Stokes simulations, implemented with the FEniCS library~\cite{alnaes2015fenics} and a nonlinear solver from the PETSc library~\cite{petsc-web-page}.
The boundary quantities are treated as variables defined over the entire domain, with values set to zero at all interior \glspl{dof}.

\subsection{2D Stokes flow}
We base the 2D Stokes flow on an analytical solution presented by Burman and Hansbo~\cite{BurmanHansbo2006}.
The exact solution is given by $\vec{v} = (20xy^3, 5x^4 - 5y^4)$ and $p = 60x^2y - 20y^3 - 5$. 
A suitable stabilization for the Stokes flow has been described in the aforementioned study~\cite{BurmanHansbo2006}, which we adopt for both 2D and 3D Stokes flow simulations in this work.
It is important to note that stabilization is necessary only when using the P1/P1 element, as the P2/P1 element is inherently stable.
The pressure term penalizes jumps in gradients of pressure, while the velocity term imposes penalty for discontinuities in divergence of velocity across elements.
It reads
\begin{gather}
    \mathcal{S}[(\vec{v},p),(\vec{w},q)] =  j(p, q) + \Tilde{j}(\vec{v}, \vec{w}), \nonumber \\
    j(p, q) = \sum_K \frac{1}{2} \int_{\partial\Omega} \gamma_p h_K^{s+1} \llbracket \vec{n}\cdot\nabla p \rrbracket \llbracket \vec{n}\cdot\nabla q \rrbracket \ds, \\
    \Tilde{j}(\vec{v}, \vec{w}) = \sum_K \frac{1}{2} \int_{\partial\Omega} \gamma_{\vec{v}} h_K^{s+1} \llbracket \nabla\cdot \vec{v} \rrbracket \llbracket \nabla\cdot \vec{w} \rrbracket \ds,
\end{gather}
where $K$ denotes a simplex in a discretized computational domain, $\llbracket x \rrbracket$ denotes the jump of quantity $x$ over an interior facet $\partial K$, $\llbracket x \rrbracket = 0$ for all exterior facets. 
The stabilization weights were set to $\gamma_p = \gamma_{\vec{v}} = 10^{-2}$.
The coefficient $s$ is defined as
\begin{equation}
    s = 
    \begin{cases}
    2 \qquad\text{if } \nu \geq h, \nonumber\\
    1 \qquad\text{if } \nu < h.
    \end{cases}
\end{equation}

Unit square meshes were generated in FEniCS using $N_i = 2^{3+i}, i=0, \dots, 6,$ elements across each side, ensuring symmetry.
When employing P2/P1 elements, we impose the velocity boundary conditions strongly.
However, when using P1/P1 with stabilization, it is important to impose the Dirichlet condition weakly to achieve optimal convergence~\cite{BurmanHansbo2014}.

Carey~\cite{Carey1985} demonstrated a significant decline in the quadratic convergence of boundary fluxes at corner elements. Specifically, the convergence rate was reduced from the expected value of $2$ to approximately $1.4$. To address this issue, our implementation of the boundary-flux evaluation method is done by computing the \gls{wss} separately on each side of the square domain. This involves subtracting the appropriate Neumann boundary integrals, following a similar approach to that presented in Equation~\eqref{eqn:wss_weak_Dirichlet} for our 3D problem. Finally, the \gls{wss} is obtained by summing the contributions from the individual boundary segments.

\subsection{3D Poiseuille flow}
We examine the Poiseuille flow through a cylindrical pipe of radius $R=1$~mm and length $L=2$~mm, which is aligned with the z-axis.
The driving force of the flow is parabolic velocity profile with a maximum value of $1$~m/s prescribed on the circular inlet. Hence, the mean velocity is $0.5$~m/s.
The no-slip \gls{bc} is imposed on the wall and do-nothing \gls{bc} together with normal outflow at the outlet
\begin{subequations}
  \begin{align}
   \hspace{4.0cm} 
    \vec{v}
    &= \vec{v_{in}} &\text{on }& \Gamma\ix{in},\\
    \vec{v}
    &= \vec{0} &\text{on }& \Gamma\ix{wall},\\
    \mathbb{T} \vec{n} 
    &= \vec{0} &\text{on }& \Gamma\ix{out},\\
    \vec v_t & = \vec{0} &\text{on }& \Gamma\ix{out}.
   \end{align}
   \label{eqn:BCs_Stokes}
\end{subequations}
The Cauchy stress tensor for an incompressible fluid flow contains the symmetric part of the velocity gradient, $\mathbb{T}(\vec v, p) = -p \mathbb{I} + 2 \nu \mathbb{D}$, rather than the complete velocity gradient, as is the case in 2D Stokes flow. 
To simulate the Poiseuille flow, it is essential to impose a boundary condition that restricts flow at the outlet to normal components only. 
For a geometry as simple as the one considered here, this is achieved by applying a Dirichlet boundary condition with zero values for the x- and y-components of velocity. 
However, for more complex geometries, a generalized approach, such as Nitsche's method, becomes necessary.
The stabilization weights for P1/P1 element were adjusted for this example and set to  $\gamma_p = 1; \gamma_{\vec{v}} = 10^{-3}$.

It can be shown that the analytical solution for \gls{wss} reads
\begin{equation}
\boldsymbol{\tau}_\text{exact} = \frac{2 \mu u_m}{R} = 8 \text{ Pa},
\end{equation}
where the maximum velocity is $u_m = 1$~m/s and the dynamic viscosity $\nu = 4 \text{ mPa}\cdot\text{s}$.

In this study we used two types of tetrahedral meshes; with and without boundary layer refinement.
Boundary layers were generated with the initial layer height set to 10\% of the surface edge length, with each subsequent layer increasing by 10\%.
A total of four boundary layers were included. To ensure consistency during mesh sensitivity analysis, the volume edge length was set equal to the surface edge length.
Meshes for studying 3D Poiseuille flow were generated using COMSOL Multiphysics 6.1 (\href{www.comsol.com}{www.comsol.com}).
Uniform meshes were used for the sensitivity analysis with the P2/P1 element, while boundary layer meshes were employed with the P1/P1 stabilized element.

For verification, we compare our results with those obtained using COMSOL Multiphysics 6.1, a commercial software. In COMSOL Multiphysics, we use the \texttt{Creep Flow} physics module, which corresponds to Stokes flow. We apply the same boundary conditions as previously described and use a direct solver to obtain the velocity and pressure fields. We then evaluate the \gls{wss} by multiplying the viscosity by the shear rate, which is automatically computed by the software.

\subsection{Patient-specific simulations}
\begin{figure}
    \centering
    \includegraphics[width=0.666\linewidth]{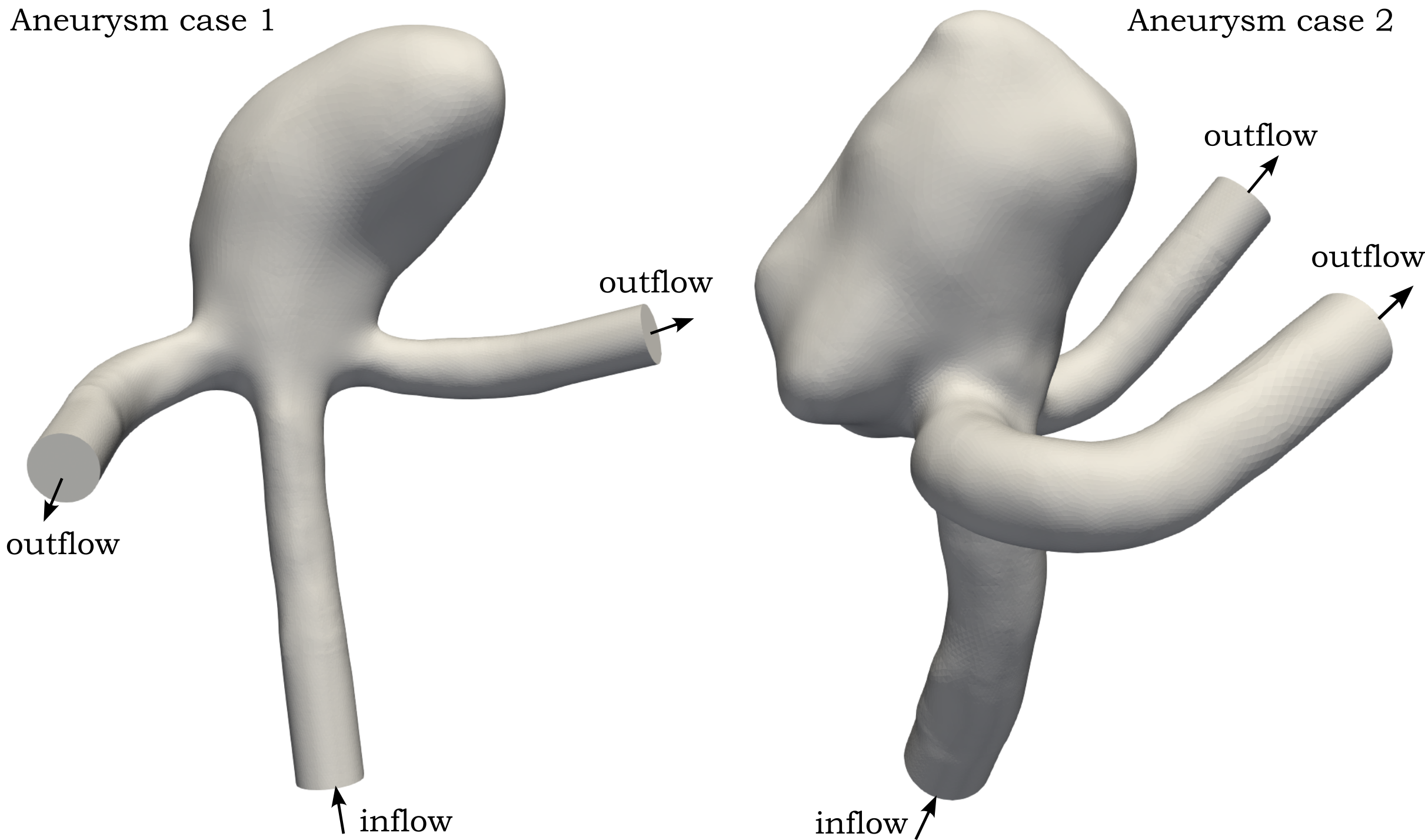}
    \caption{Computational geometries for both aneurysm cases considered in this study, with the inflow branch and two outflow branches indicated by arrows.}
    \label{fig:meshes}
\end{figure}

Two patient-specific aneurysm geometries were adopted from~\cite{Janiga2015}; both aneurysm geometries can be seen in ~\Cref{fig:meshes}.
Surface remeshing for these geometries was performed using the PyMeshlab~\cite{pymeshlab} package, achieving the target surface edge length, which ranged from $0.3$~mm for the coarsest mesh to $0.1$~mm for the finest mesh. Cylindrical flow extensions were added to both inflow and outflow branches.
Subsequently, the open-source GMSH library~\cite{gmsh2009} is employed to generate volume meshes consistent with the previous example, producing one mesh type with boundary layer refinement and another without.
The number of boundary layers is again set to four, with the first layer height at 10\% of the surface edge length and each subsequent layer increasing by 10\%.
Meshes with boundary layers were employed exclusively for the P1/P1 stabilized element, whereas uniform meshes were employed for the P2/P1 element to enable a consistent mesh sensitivity test, which would not have been feasible with boundary layer meshes due to computational costs associated with monolithic solver.

Although the numerical solver is inherently designed for time-dependent simulations, a steady-state solution can be found using a suitable adaptive time-stepping.
In particular, we use the BDF2 time-stepping with adaptivity in time steps between $10^{-2} - 10^5$. 
Initially, the velocity field is set to zero throughout the domain.
Subsequently, the inflow velocity is linearly ramped up over a period of $0.5$ seconds to initiate the desired inflow, which is parabolic velocity profile with a mean value of $0.5$~m/s, consistent with~\cite{valen2013high}. The corresponding Reynolds numbers are $269$ for aneurysm case 1 and $263$ for aneurysm case 2.
The total simulation length is set to $T =10^6$ seconds. These settings ensure that the solution effectively reaches a steady state.

For P1/P1 element, we use the interior penalty stabilization as introduced in~\cite{Burman2007}
\begin{gather}
\mathcal{S}[(\vec{v},p),(\vec{w},q)] = \nonumber \\
    \sum_{K \in \mathcal{T}_h} \left[ \alpha_i h^2 \int_{\partial K} |\vec v \cdot \vec n|^2 \llbracket\nabla \vec v\rrbracket \llbracket\nabla \vec w\rrbracket \text{dS}
    + \alpha_{\vec{v}} h^2 \int_{\partial K} \llbracket\nabla \vec v\rrbracket \llbracket\nabla \vec w\rrbracket \text{dS}
    + \alpha_p h^2 \int_{\partial K} \llbracket\nabla p\rrbracket \llbracket\nabla q\rrbracket \text{dS} \right],
\end{gather}
where the weights were chosen in accordance with the Poiseuille flow problem above, $\alpha_{\vec{v}} = 10^{-3}, \alpha_p = 1$, and $\alpha_i = 10^{-3}$.

For each finite element selection, the mean \gls{wss} across the aneurysm dome is assessed using both the projection method and the boundary-flux evaluation technique, as described in Sections \ref{sec:projection_method} and \ref{sec:boundary-flux_evaluation}, respectively. 
Again, we consider the \gls{dg0}, \gls{dg1} and \gls{cg1} spaces for the projection method. 
Subsequently, for each \gls{wss} evaluation method we compute the maximum, minimum and average \gls{wss} values over the aneurysm dome as well as the \gls{lsa} indicator. The \gls{lsa} is defined as the percentage of the aneurysm wall area exposed to \gls{wss} below $10\%$ of the mean \gls{wss} in the parent artery, as introduced in~\cite{Jou2008}.

Post-processing computations are performed using the FEniCS library\cite{alnaes2015fenics}, interactive mesh clipping and surface plots are conducted in ParaView (\href{https://www.paraview.org/}{www.paraview.org}).

\section{Results}\label{results}

\subsection{2D Stokes flow}
\begin{figure}[htbp]
    \centering
    \includegraphics[width=\linewidth]{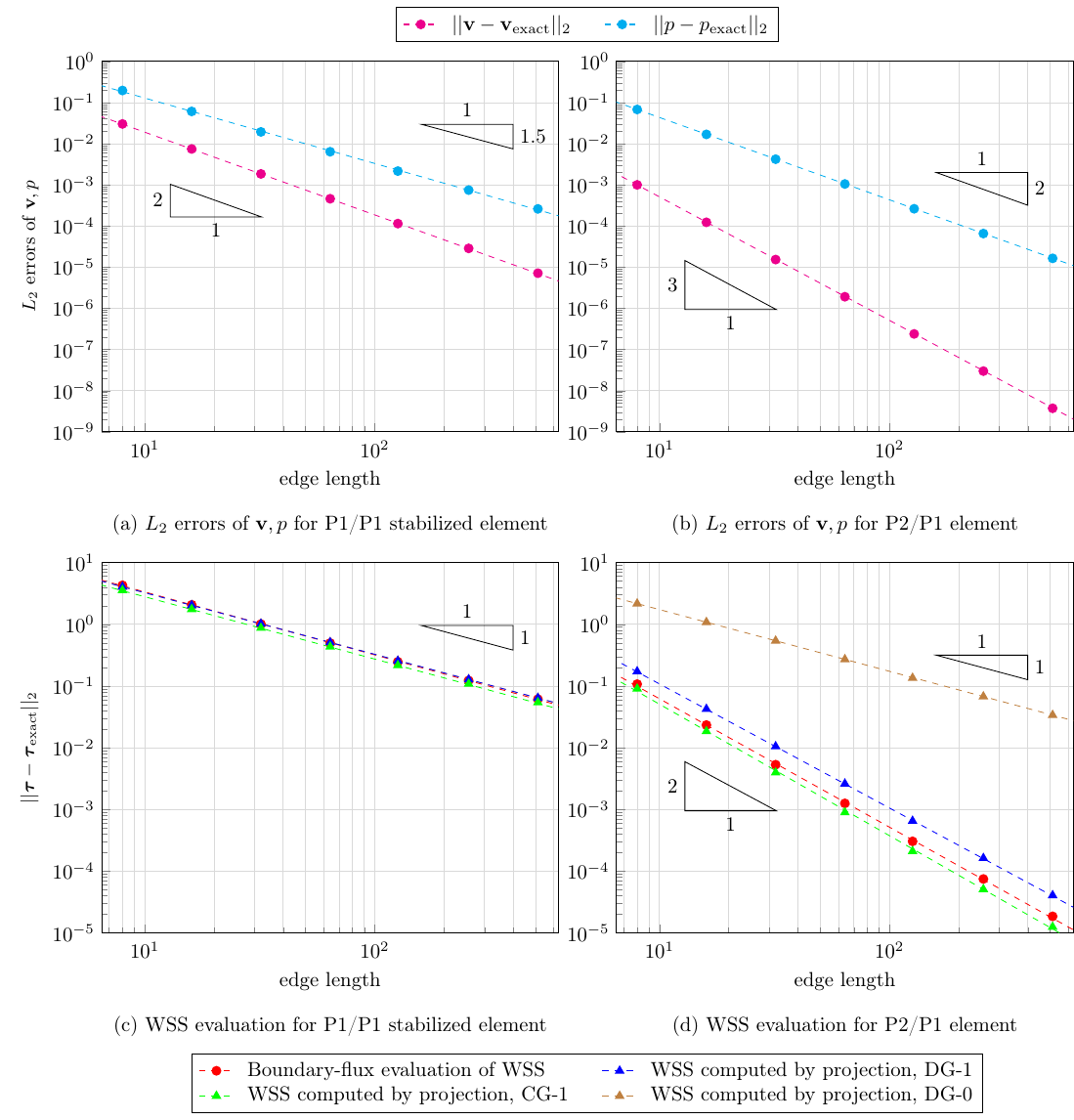}
    \caption{2D Stokes flow: Logarithmic plot of $L_2$ error in velocity, pressure, and \gls{wss} as a function of edge length, with errors evaluated against the analytical solution. Red circles show boundary-flux evaluation; green, blue and brown triangles show \gls{wss} computed by projection in \gls{cg1}, \gls{dg1} and \gls{dg0} space, respectively. Dashed lines represent computed convergence rate for corresponding datasets. 
    (a) $L_2$ errors of $\vec v, p$ for P1/P1 stabilized element.
    (b) $L_2$ errors of $\vec v, p$ for P2/P1 element.
    (c) \gls{wss} assessment for P1/P1 stabilized element.
    (d) \gls{wss} assessment for P2/P1 element.}
    \label{fig:2d_stokes}
\end{figure}

\Cref{fig:2d_stokes} shows $L_2$ errors of velocity, pressure, together with all considered \gls{wss} assessment methods. Errors were calculated with respect to analytical solutions.
When using P1/P1 stabilized element, the convergence rate is $2$ for velocity and $1.59$ for pressure, see \Cref{fig:2d_stokes}(a).
While for P2/P1 element, convergence rate of $3$ is achieved for velocity and $2$ for pressure as shown in \Cref{fig:2d_stokes}(b).
\gls{wss} assessment for P1/P1 stabilized element yielded linear convergence and the differences between the methods are very small.
For P2/P1 element, the convergence is quadratic for both boundary-flux evaluation and standard evaluation (\gls{dg1} or \gls{cg1}).
However, only linear convergence is achieved for \gls{dg0} \gls{wss}.  

\subsection{3D Poiseuille flow}

\begin{figure}[htbp]
    \centering
    \includegraphics[width=\linewidth]{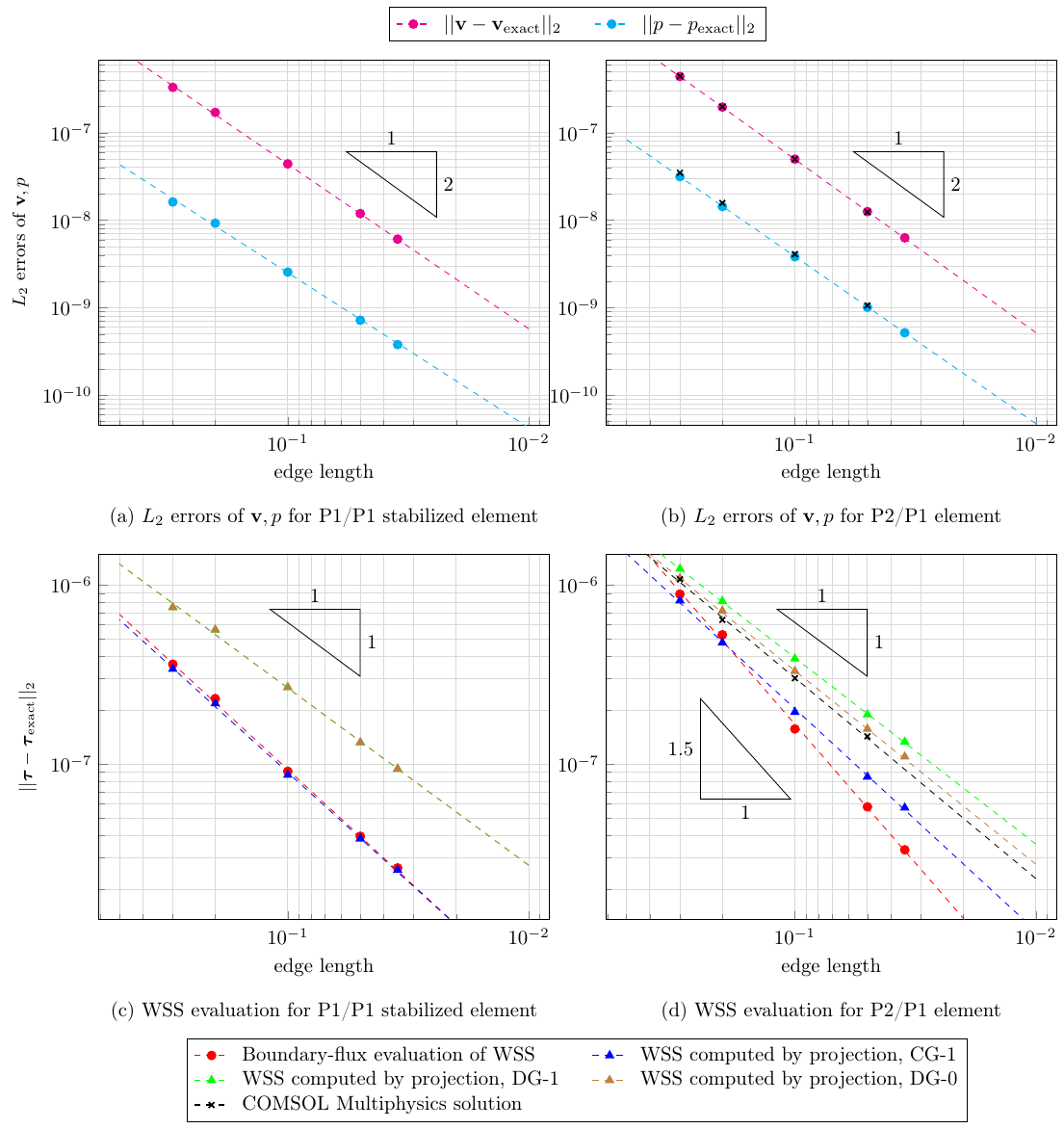}
    \caption{3D Poiseuille flow: Logarithmic plot of $L_2$ error in velocity, pressure, and \gls{wss} as a function of edge length, with errors evaluated against the analytical solution. Red circles show boundary-flux evaluation; green, blue and brown triangles show \gls{wss} computed by projection in \gls{cg1}, \gls{dg1} and \gls{dg0} space, respectively. Dashed lines represent computed convergence rate for corresponding datasets. 
    (a) $L_2$ errors of $\vec v, p$ for P1/P1 stabilized element on meshes with boundary layers.
    (b) $L_2$ errors of $\vec v, p$ for P2/P1 element on uniform meshes.
    (c) \gls{wss} assessment for P1/P1 stabilized element on meshes with boundary layers.
    (d) \gls{wss} assessment for P2/P1 element on uniform meshes.
    }
    \label{fig:loglog_plot_poiseuille}
\end{figure}

Figure \ref{fig:loglog_plot_poiseuille} shows $L_2$ errors of velocity, pressure and \gls{wss} with respect to analytical solutions.
The results for the P1/P1 stabilized element, which employ meshes with boundary layers, are shown on the left.
The results for the P2/P1 element, which uses uniform meshes, are displayed on the right.
Since the momentum equation is divided by density, as is typical for Stokes problems, it should be noted that both the pressure and \gls{wss} values presented in the Poiseuille flow example are scaled by density.

The convergence rate for velocity is $1.88$ and $1.78$ for pressure in case of P1/P1 stabilized element. Using P2/P1 element and our open-source FEniCS code, velocity and pressure converge with rates of $1.98$ and $1.91$, respectively, and COMSOL Multiphysics solutions give very similar rates of $2.00$ and $1.95$, respectively.

\gls{wss} evaluation for the P1/P1 stabilized element yields a rate of $1.24$ for the boundary-flux evaluation, $1.22$ for the projection into \gls{cg1} space, and $0.99$ for both projections to \gls{dg0} and \gls{dg1}.
For the P2/P1 element, a significant improvement in convergence rate is observed for the boundary-flux evaluation, reaching a rate of $1.55$. The standard evaluation attains a rate of $1.24$ using \gls{cg1} space, while $1.04$ and $1.08$ for \gls{dg1} and \gls{dg0} spaces. \gls{wss} evaluated in COMSOL converges with a rate of $1.12$.
All $L_2$ errors of velocity, pressure and \gls{wss}, together with corresponding convergence rates, can be found in Tables S1-S2 in the Supplementary materials.

\subsection{Patient-specific simulations}
Maximum, minimum and average \gls{wss} over the aneurysm dome were assessed for each choice of finite element and for both the projection and the boundary-flux evaluation method; see~\Cref{fig:WSS_case01,fig:WSS_case02}.
Additionally, all \gls{wss} values are shown in in Tables S3-S5 in the Supplementary materials, together with relative percentage differences with respect to the \gls{cg1} projection method for each mesh size, as this method is most frequently used.
Again, the P2/P1 results were computed on uniform meshes, while P1/P1 stabilized element was used for meshes with boundary layers.
It can be observed that the P2/P1 element is much more robust in assessing the average values of \gls{wss} than the P1/P1 stabilized element.

\begin{figure}[htbp]
    \centering
    \includegraphics[width=.9\linewidth]{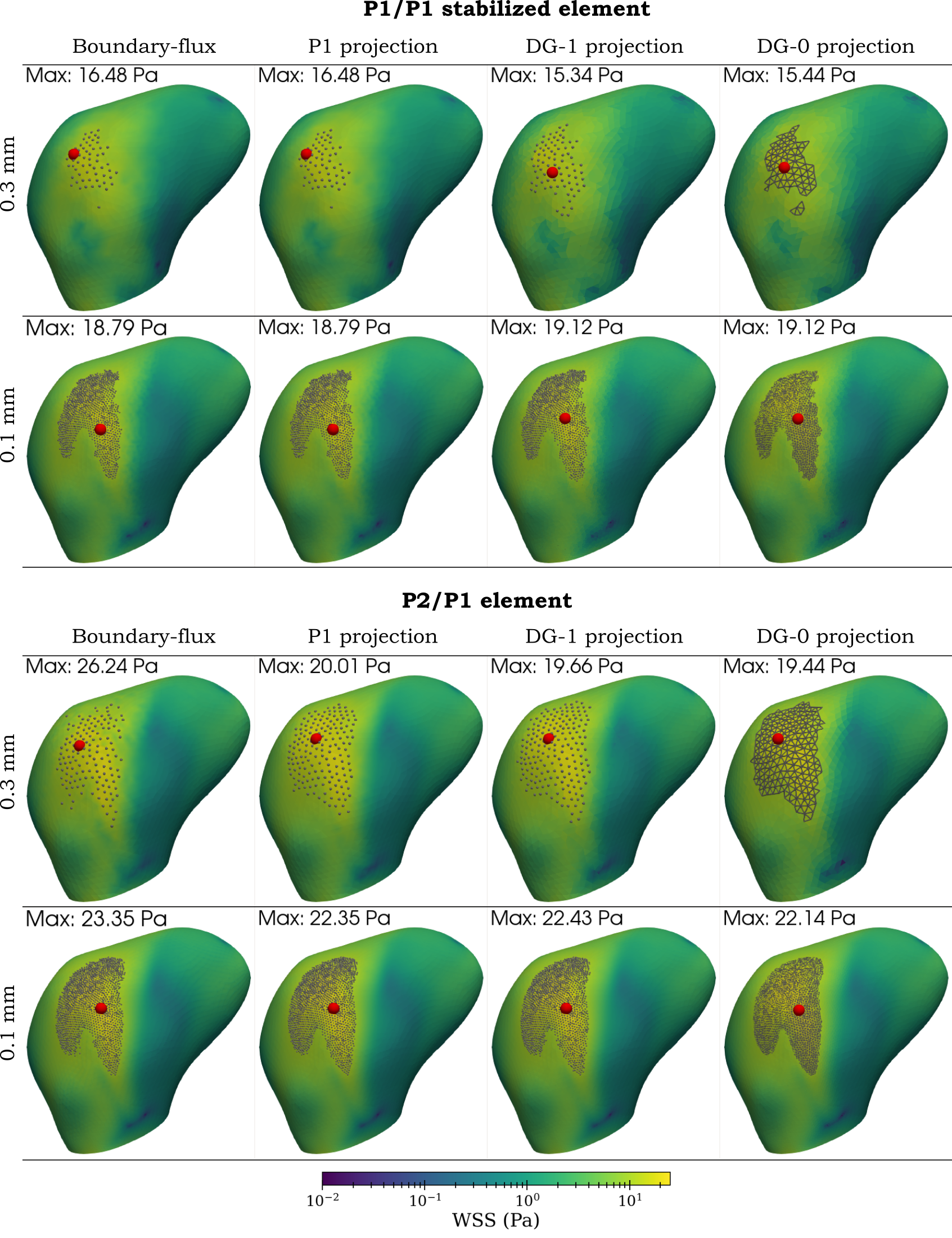}
    \caption{Aneurysm case 1: Comparison of boundary-flux evaluation and projection method for P2/P1 element. Red points denote the maximum value of \gls{wss} over the aneurysm dome, and gray points or triangles represent vertices or facets where $\text{WSS}\geq10$~Pa. \textit{Top row}: edge length $0.3$~mm; \textit{bottom row}: edge length $0.1$~mm.}
    \label{fig:wss_case01_paraview}
\end{figure}

\begin{figure}[htbp]
    \centering
    \includegraphics[width=.9\linewidth]{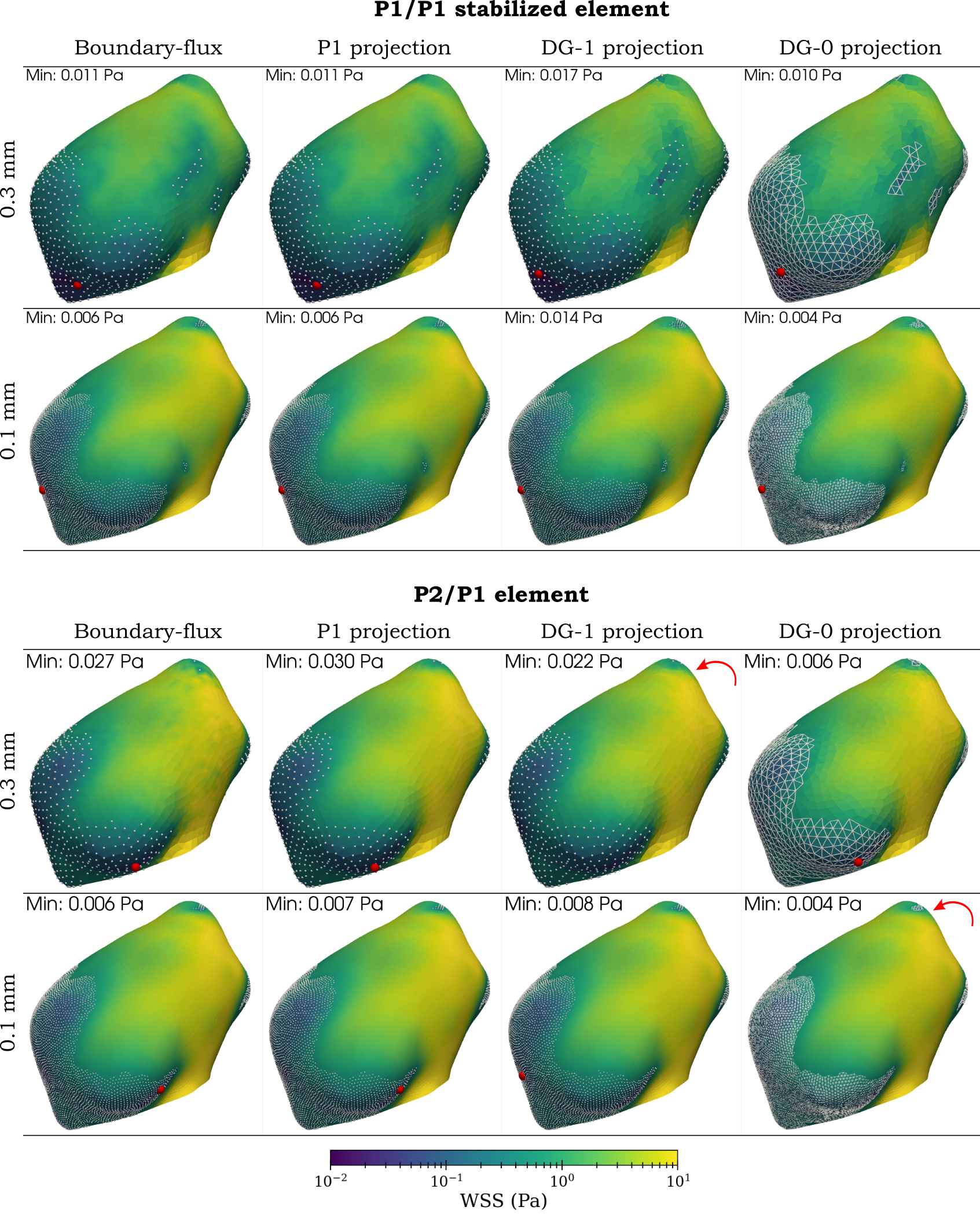}
    \caption{Aneurysm case 2: Comparison of boundary-flux evaluation and projection method. Red points denote the minimum value of \gls{wss} over the aneurysm dome, red arrows indicate the minimum is located outside the current viewpoint; and white points or triangles represent vertices or facets where $\text{WSS}\leq 0.5$~Pa. \textit{Top row}: edge length $0.3$~mm; \textit{bottom row}: edge length $0.1$~mm.}
    \label{fig:wss_case02_paraview}
\end{figure}

Moreover, a visual comparison of all \gls{wss} assessment techniques is depicted in Figures~\ref{fig:wss_case01_paraview}-\ref{fig:wss_case02_paraview}.
The first of these figures shows the case 1 aneurysm with an annotated maximum \gls{wss} over the dome and all boundary vertices (or facets) where \gls{wss} is above $10$~Pa. The second figure shows the case 2 aneurysm with an annotated minimum value over the dome and all boundary vertices (or facets) where \gls{wss} is below $0.5$~Pa.

\begin{figure}[htbp]
    \centering
    \includegraphics{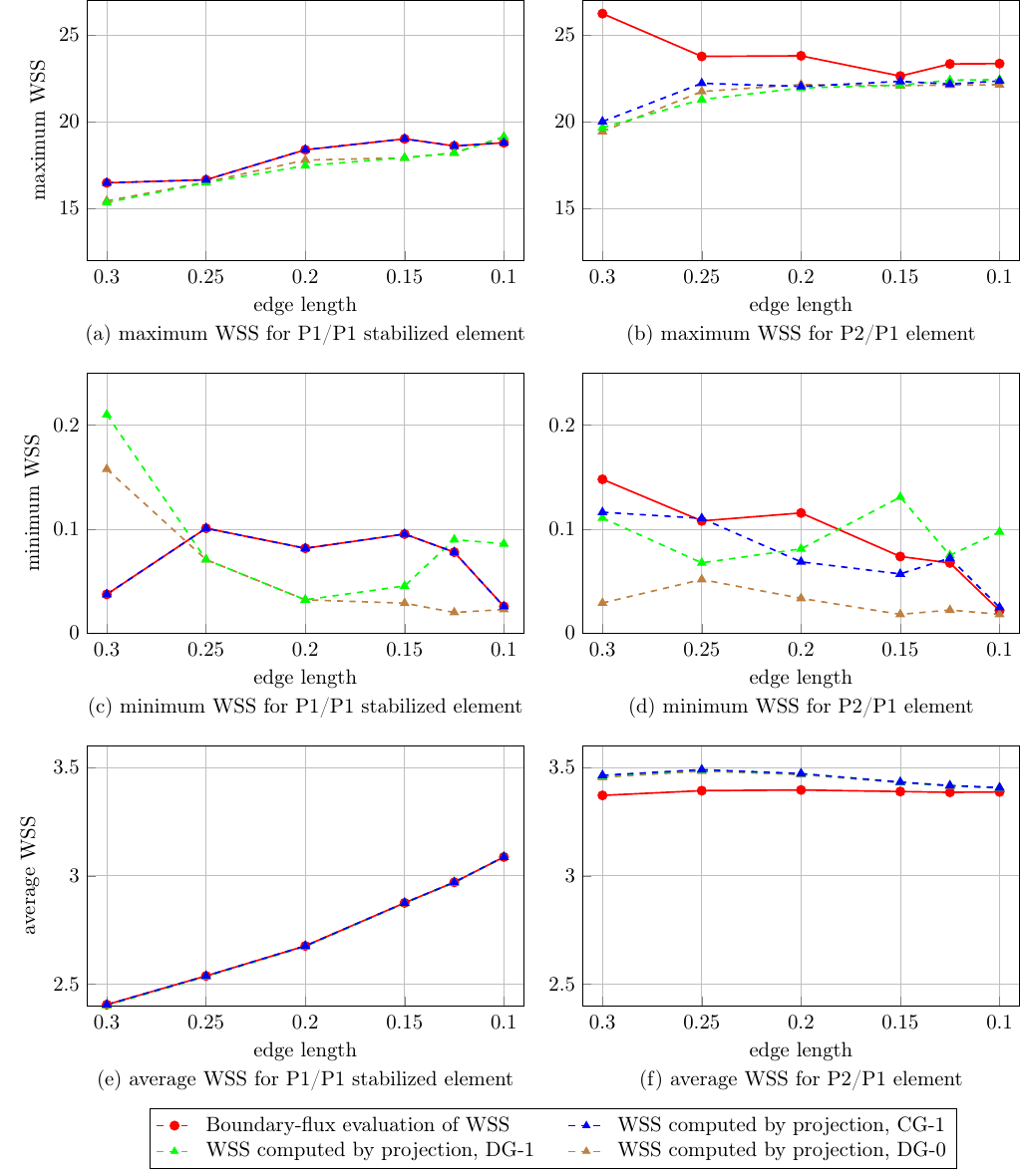}
    \caption{Aneurysm case 1: Maximum, minimum and average values of \gls{wss} over the aneurysm dome for each \gls{wss} evaluation method. (a, c, e) P1/P1 stabilized element, evaluated on meshes with boundary layers. (b, d, f) P2/P1 element, evaluated on uniform meshes.}
    \label{fig:WSS_case01}
\end{figure}

\begin{figure}[htbp]
    \centering
    \includegraphics{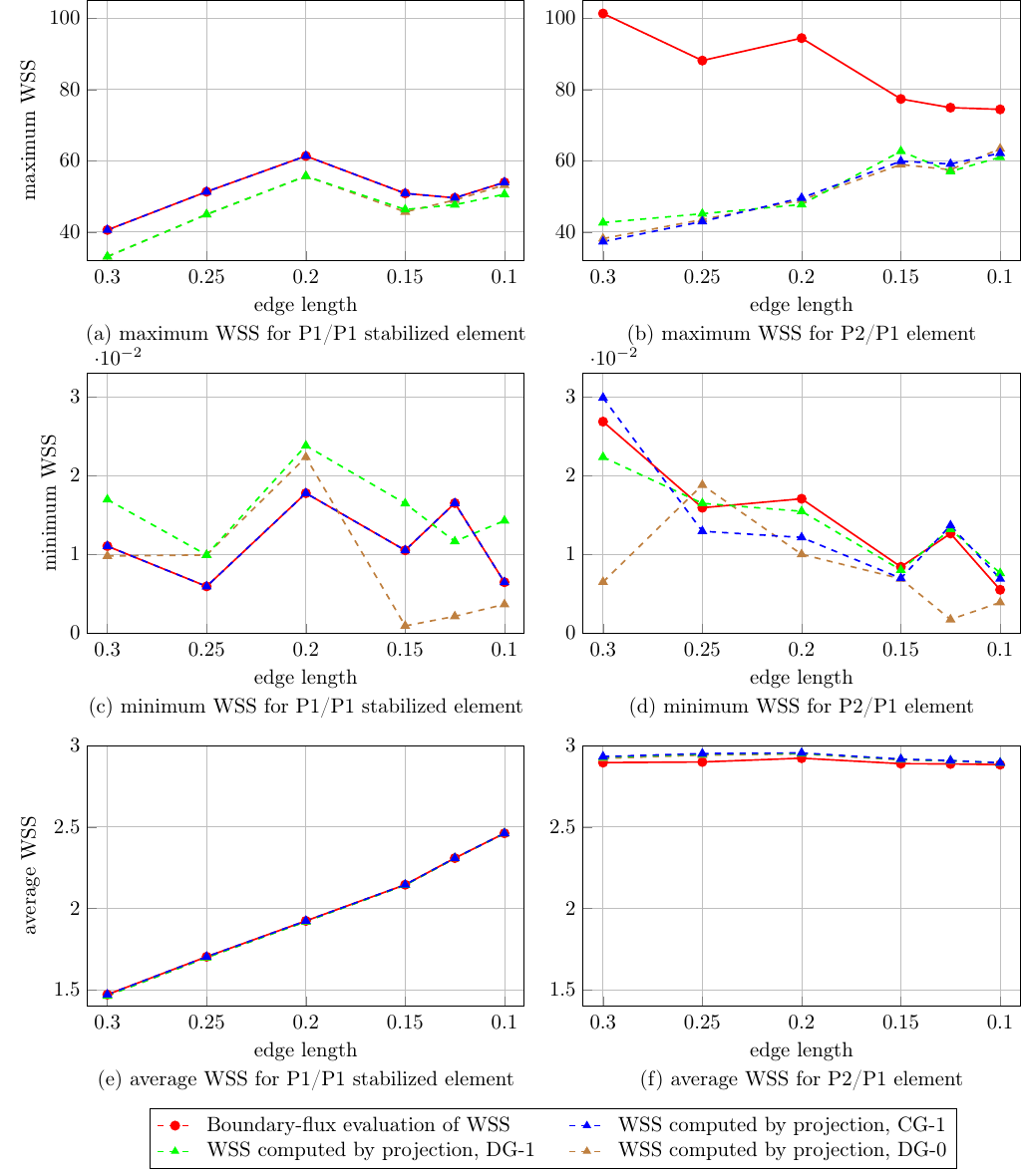}
    \caption{Aneurysm case 2: Maximum, minimum and average values of \gls{wss} over the aneurysm dome for each \gls{wss} evaluation method. (a, c, e) P1/P1 stabilized element, evaluated on meshes with boundary layers. (b, d, f) P2/P1 element, evaluated on uniform meshes.}
    \label{fig:WSS_case02}
\end{figure}

\Cref{tab:LSA_case01} presents the values of \gls{lsa} obtained through simulations using either the P1/P1 stabilized element or the P2/P1 element. The \gls{wss} is assessed using all considered evaluation techniques.
When using the P1/P1 stabilized element, the \gls{lsa} indicator exhibits a strong dependence on mesh size, with differences exceeding $100\%$ between the coarsest and finest meshes in case 1 aneurysm. In contrast, the variation between individual evaluation methods is relatively minor (less than $1\%$). Conversely, for the P2/P1 element, the \gls{lsa} indicator demonstrates robustness to changes in mesh size.
However, in this case, the boundary-flux evaluation deviates more significantly from the projection methods, with differences up to $5\%$.

\definecolor{mygray}{gray}{0.9}

\begin{center}
\begin{table*}[!h]%
\caption{Hemodynamic indicator \gls{lsa} $[\%]$ for P1/P1 stabilized element and P2/P1 element across all \gls{wss} evaluation methods. \label{tab:LSA_case01}}
\begin{tabular*}{\textwidth}{@{\extracolsep\fill}lcccccccc@{}}
\toprule
\multicolumn{9}{@{}c}{\textbf{P1/P1 stabilized element; meshes containing boundary layers}} \\ 
\midrule
\multirow{2}{*}{\textbf{edge length [mm]}}  & \multicolumn{2}{@{}c}{\textbf{boundary-flux evaluation}}  & \multicolumn{2}{@{}c}{\textbf{\gls{cg1} projection}}  & \multicolumn{2}{@{}c}{\textbf{\gls{dg1} projection}}  &  \multicolumn{2}{@{}c}{\textbf{\gls{dg0} projection}} \\
 & case 1 & case 2 & case 1 & case 2 & case 1 & case 2 & case 1 & case 2 \\
\hline
0.300 & \cellcolor{mygray} 16.86 & 47.83 & \cellcolor{mygray} 16.86 & 47.83 & \cellcolor{mygray} 16.76 & 47.85 & \cellcolor{mygray} 16.76 & 47.85 \\
0.250 & \cellcolor{mygray} 11.37 & 42.46 & \cellcolor{mygray} 11.37 & 42.46 & \cellcolor{mygray} 11.46 & 42.40 & \cellcolor{mygray} 11.46 & 42.40 \\
0.200 & \cellcolor{mygray} 9.55 & 37.96 & \cellcolor{mygray} 9.55 & 37.96 & \cellcolor{mygray} 9.46 & 38.04 & \cellcolor{mygray} 9.46 & 38.03 \\
0.150 & \cellcolor{mygray} 8.38 & 36.05 & \cellcolor{mygray} 8.40 & 36.05 & \cellcolor{mygray} 8.31 & 36.10 & \cellcolor{mygray} 8.31 & 36.10 \\
0.100 & \cellcolor{mygray} 8.14 & 35.00 & \cellcolor{mygray} 8.14 & 34.99 & \cellcolor{mygray} 8.12 & 35.09 & \cellcolor{mygray} 8.12 & 35.09 \\
\hline
\midrule
\multicolumn{9}{@{}c}{\textbf{P2/P1 element; uniform meshes}} \\
\midrule
\multirow{2}{*}{\textbf{edge length [mm]}}  & \multicolumn{2}{@{}c}{\textbf{boundary-flux evaluation}}  & \multicolumn{2}{@{}c}{\textbf{\gls{cg1} projection}}  & \multicolumn{2}{@{}c}{\textbf{\gls{dg1} projection}}  &  \multicolumn{2}{@{}c}{\textbf{\gls{dg0} projection}} \\
 & case 1 & case 2 & case 1 & case 2 & case 1 & case 2 & case 1 & case 2 \\
\hline
0.300 & \cellcolor{mygray} 7.82 & 34.81 & \cellcolor{mygray} 7.53 & 33.38 & \cellcolor{mygray} 7.49 & 33.35 & \cellcolor{mygray} 7.55 & 33.21 \\
0.250 & \cellcolor{mygray} 7.77 & 34.79 & \cellcolor{mygray} 7.75 & 33.92 & \cellcolor{mygray} 7.66 & 33.84 & \cellcolor{mygray} 7.48 & 33.74 \\
0.200 & \cellcolor{mygray} 7.80 & 34.58 & \cellcolor{mygray} 7.70 & 34.33 & \cellcolor{mygray} 7.79 & 34.41 & \cellcolor{mygray} 7.64 & 34.39 \\
0.150 & \cellcolor{mygray} 7.80 & 34.42 & \cellcolor{mygray} 7.70 & 34.27 & \cellcolor{mygray} 7.67 & 34.28 & \cellcolor{mygray} 7.66 & 34.26 \\
0.100 & \cellcolor{mygray} 7.78 & 34.49 & \cellcolor{mygray} 7.68 & 34.38 & \cellcolor{mygray} 7.68 & 34.41 & \cellcolor{mygray} 7.68 & 34.39 \\
\hline
\bottomrule
\end{tabular*}
\end{table*}
\end{center}

\section{Discussion}\label{discussion}
\subsection{2D Stokes flow}
As a proof of concept for the methods presented in this paper, we investigated a 2D Stokes flow on a unit square.
The analytical solution is straightforward and it exhibits either constant or linearly increasing behavior along the sides, with a discontinuity located at the upper right corner.
Although theoretical predictions indicate quadratic convergence for both velocity and pressure when using the P1/P1 stabilized element \cite{BurmanHansbo2006}, tuning the stabilization parameters to achieve this level of convergence proved to be quite challenging. In our experiments, we achieved quadratic convergence for velocity but only a convergence rate of $1.59$ for pressure.
In contrast, employing the P2/P1 element resulted in convergence rates of $3$ for velocity and $2$ for pressure, which is in accordance with the theoretical results.

\gls{wss} assessment for this problem demonstrated linear convergence for all methods when using the P1/P1 stabilized element, with no significant differences observed in the absolute errors among them.
In contrast, when employing the P2/P1 element, linear convergence was achieved only with the \gls{dg0} projection.
For the \gls{dg1} and \gls{cg1} projections, as well as boundary-flux evaluation, quadratic convergence is attained.

It is important to note that convergence rates can be significantly influenced by the nature of the analytical solution.
In this example, the analytical \gls{wss} on the upper boundary exhibits a linearly increasing behavior, which likely explains why the convergence rate of the \gls{dg0} projection is one order lower than that of the other methods.

\subsection{3D Poiseuille flow}

Using the P1/P1 stabilized element, similar convergence rates are observed for the boundary-flux evaluation ($1.24$) and \gls{cg1} projection method ($1.22$), while the \gls{dg1} and \gls{dg0} projection converge with a rate of $0.99$.    
For the P2/P1 element, the difference between the two methods becomes more prominent. 
The boundary-flux evaluation achieves a convergence rate of $1.55$, compared to $1.24, 1.04$ and $1.08$ for the standard \gls{cg1}, \gls{dg1} and \gls{dg0} projections.

Similarly to the previous example, the nature of the analytical solution, which in this case is a constant function, may affect the convergence rates of the studied methods. 
Additionally, the convergence is influenced by the approximation of the cylindrical computational domain using only linear tetrahedral elements.
A possible improvement of \gls{wss} assessment may include using higher order boundary elements.
Although it would be straightforward to do such improvement of the boundary of a cylinder, there is no simple extension to patient-specific geometries, which would limit the applicability of such method in practice.

\subsection{Patient-specific simulations}

Using the P1/P1 stabilized element, the boundary-flux evaluation and the \gls{cg1} projection method yield equivalent results, which is in accordance with the Poiseuille flow example.
Discontinuous projections yield differences of maximum values with respect to the continuous projection up to $6.9\%$ (case 1) or $18.5\%$ (case 2) on the coarsest mesh and $1.8\%$ (case 1) or $-6.9\%$ (case 2) on finest mesh.
Using the P2/P1 element, the boundary-flux evaluation differs from the projection methods, particularly in determining maximum values.
The boundary-flux method seems to be more susceptible to inaccuracies in the pressure field, which propagate to the assessment of \gls{wss}, unlike the projection method.
Introducing a pressure stabilization term might be beneficial to improve the boundary-flux evaluation, especially on coarse meshes. Alternatively, employing boundary layer meshes for P2/P1 elements might be considered, although this approach would significantly increase computational costs. However, as expected, these discrepancies diminish with mesh refinement.

Furthermore, the results indicate that the minimum \gls{wss} value over the aneurysm dome may not be a reliable hemodynamic indicator. Its location shows high sensitivity to both the mesh size and the mixed finite element space employed in the simulation (see \Cref{fig:wss_case02_paraview}). This sensitivity affects the consistency and reproducibility of using the minimum \gls{wss} value as a hemodynamic indicator in clinical or research applications.

A commonly used hemodynamic indicator derived from \gls{wss} is the \gls{lsa}, which captures the percentage of the aneurysm dome exposed to low wall shear stress (below $10\%$ of the mean value in the parent artery).
Our results suggest that this indicator is robust to mesh size when P2/P1 element is used.
However, over $100\%$ differences in \gls{lsa} were observed for P1/P1 stabilized element in aneurysm case 1, which is related to the observation that the average \gls{wss} is largely underestimated for coarse meshes, see~\Cref{fig:WSS_case01}(e) and \Cref{fig:WSS_case02}(e).

\subsection{Limitations}
As highlighted in the review paper by Steinman et al. \cite{Steinman2019}, numerous factors influence the patient-specific modeling pipeline, many of which are now well understood and should be standardised. In this context, the most notable limitation of our study is the use of only two anatomically realistic geometries, which precludes any statistical analysis. Additionally, we focused exclusively on steady-state flows, limiting the flow regime to laminar, although transitional flows could be present in 50\% of all aneurysms~\cite{khan2021prevalence}, if care is taken to resolve them~\cite{valen2014mind}. Consequently, we do not claim that the simulations are patient-specific, although steady-state simulations have been shown to approximate time-averaged WSS~\cite{geers2014approximating,cebral2011quantitative}.

\subsection{Implications and future work}

From a mathematical perspective, discontinuous spaces for \gls{wss} may appear more natural, given the lack of guaranteed continuity in velocity gradients.
Although physiological contexts suggest that stresses should ideally be continuous, projecting discontinuous finite element results into continuous spaces, such as \gls{cg1}, can introduce artifacts like the Gibbs phenomenon, as demonstrated in \cite{Zhang2022}. Therefore, it warrants attention whether the choice of finite element formulation introduces such artifacts, as this may impact the reliability of the results.

A more general boundary condition than the standard no-slip condition can be considered for the walls, such as Navier's slip boundary condition. 
This condition establishes a relationship between the tangential velocity on the boundary and the tangential traction (i.e., \gls{wss}). 
With such a boundary condition, \gls{wss} assessment becomes particularly straightforward, as \gls{wss} is essentially a rescaled tangential velocity at the boundary, see \cite{Chabiniok2021}. 
Consequently, there is no need to solve an additional system to evaluate \gls{wss}. 
From this perspective, \gls{wss} could live in the same finite element space as the velocity field.

In this work, we used the \gls{cg1} space for the boundary-flux evaluation of \gls{wss} with both the P1/P1 stabilized element and the P2/P1 element. The P2/P1 element results could potentially be extended to P2 boundary-flux assessments, however, it is beyond the scope of this study.

\section{Conclusions}\label{conclusions}
\gls{wss} is one of the most commonly studied hemodynamic indicators, although the methods for its numerical evaluation has not yet been studied in detail.
This study investigated numerical approaches for evaluating \gls{wss} using \gls{fem}. 
Two benchmark problems and two patient-specific aneurysm geometries were examined, focusing on P1/P1 stabilized and P2/P1 mixed finite elements for velocity and pressure. 
\gls{wss} was evaluated using a boundary-flux method introduced in this paper and the standard $L_2$ projection method applied to three finite element spaces: \gls{cg1}, \gls{dg1}, and \gls{dg0}. 

For the P1/P1 stabilized element, the boundary-flux method and the \gls{cg1} projection method yielded equivalent results. 
In contrast, with the P2/P1 element, the boundary-flux method achieved a faster convergence rate for \gls{wss} in the Poiseuille flow example; however, when applied to patient-specific geometries, it was more sensitive to inaccuracies in the pressure field compared to the projection method.

This study represents an initial step toward standardizing \gls{wss} assessment in patient-specific geometries, which could lead to more reliable computational studies applicable to clinical practice. 
Additionally, it aims to encourage discussion within the scientific community regarding the need for standardized and reliable evaluation of hemodynamic indicators.

\section*{Acknowledgments}
Jana Brunátová and Jaroslav Hron were supported by the Czech Health Research Council No. NU22-08-00124. In addition, Jana Brunátová was supported by the Grant Agency of Charles University (project no. 308522) and the project SVV-2023-260711.
J\o{}rgen S. Dokken and Kristian Valen-Sendstad were supported by Simula Research Laboratory.
Our sincere thanks go to Erik Burman for providing his insight necessary to advance this research project.

\section*{Data availability statement}
The source code for this study is publicly accessible on Zenodo \href{https://doi.org/10.5281/zenodo.14506052}{doi.org/10.5281/zenodo.14506052} and all computational meshes are available at \href{https://doi.org/10.5281/zenodo.14503385}{doi.org/10.5281/zenodo.14503385}.

\printbibliography






\end{document}